\documentclass[11pt,oneside]{article}
\usepackage[OT4]{fontenc}\usepackage[cp1250]{inputenc}\usepackage{a4} \usepackage{amsfonts}\usepackage{amsmath}\usepackage{amssymb}\usepackage{amsxtra}
\usepackage{amsthm}

\usepackage{a4}

\usepackage[dvips]{graphicx}\usepackage{pstricks}\usepackage{color}
\usepackage{bezier}\usepackage{epic}\usepackage{eepic}
\usepackage{pstricks}\usepackage{pst-node}\usepackage{lineno}
\usepackage[normalem]{ulem}\usepackage{epsf}

\newtheorem{thm}{Theorem}[section]

\newtheorem{prop}[thm]{Proposition}
\newtheorem{obs}[thm]{Observation}
\newtheorem{cor}[thm]{Corollary}

\newcommand{\cp}{\, \square \,}

\let\oldenumerate\enumerate
\renewcommand{\enumerate}{
  \oldenumerate
  \setlength{\itemsep}{0.5pt}
  \setlength{\parskip}{0pt}
  \setlength{\parsep}{0pt}
}

\begin{document}

\title{On $\alpha$-excellent graphs}
\author{$^1$Magda Dettlaff, $^{2}$Michael A. Henning\thanks{Research
supported in part by the University of Johannesburg} \, and \, $^{3}$Jerzy Topp\\
\\
$^1$Faculty of Mathematics, Physics and Informatics\\
University of Gda\'nsk \\
80-952 Gda\'nsk, Poland \\
\small \tt Email:  magda.dettlaff@pg.edu.pl \\
\\
$^2$Department of Mathematics and Applied Mathematics\\ University of Johannesburg \\ Auckland Park 2006, South Africa \\ \small {\tt mahenning@uj.ac.za} \\
\\
$^3$Institute of Applied Informatics\\ University of Applied Sciences\\ 82-300 Elbl\c{a}g, Poland \\ \small {\tt j.topp@ans-elblag.pl}}

\date{}
\maketitle

\begin{abstract} \noindent
A graph $G$ is $\alpha$-excellent if every vertex of $G$ is contained in some maximum independent set of $G$. In this paper, we characterize $\alpha$-excellent bipartite graphs, $\alpha$-excellent unicyclic graphs,  $\alpha$-excellent simplicial graphs,  $\alpha$-excellent chordal graphs,  $\alpha$-excellent block graphs, and we show that every generalized Petersen graph is $\alpha$-excellent.
\end{abstract}


\section{Introduction}

A \emph{dominating set} of a graph $G$ is a set $D$ of vertices of $G$ such that every vertex belonging to $V(G)\setminus D$ is adjacent to at least one vertex in $D$. A subset $I$ of $V(G)$ is \emph{independent} if no two vertices belonging to $I$ are adjacent in $G$. An \emph{independent dominating set}, abbreviated ID-set, of $G$ is a set that is both a dominating set and an independent set. The domination number $\gamma(G)$ of $G$ is the minimum cardinality of a dominating set in $G$, while the independent domination number $i(G)$ of $G$ is the minimum cardinality of an ID-set in $G$. The \emph{independence number} $\alpha(G)$ of $G$ is the maximum cardinality of an independent set in $G$. An $\alpha$-\emph{set of $G$} is an independent set of $G$ of maximum cardinality~$\alpha(G)$. The study of independent sets in graphs was begun by Berge~\cite{Berge1962,Berge1981} and Ore~\cite{Ore1962}. We refer the reader to the survey~\cite{GoHe-13} of results on independent domination in graphs by Goddard and Henning. For recent books on domination in graphs, we refer the reader to~\cite{HaHeHe-20,HaHeHe-21,HaHeHe-22}.

Recently another parameter concerning the existence and cardinality of independent sets in a~graph, the common independence number of a graph, was introduced and studied by Dettlaff et al.~\cite{DeLeTo-21}.  Formally, the \emph{common independence number} of a graph $G$, denoted by $\alpha_c(G)$, is defined as the greatest integer $r$ such that every vertex of $G$ belongs to some independent subset $X$ of $V(G)$ with $|X| \ge r$. Thus, the common independence number of $G$ refers to numbers of mutually independent vertices of $G$ and it emphasizes the notion of the individual independence of a vertex of $G$ from other vertices of $G$. It follows immediately from the above definitions that the common independence number is sandwiched between the independent domination number and the independence number. We state this formally as follows.

\begin{obs}{\rm (\cite{DeLeTo-21})}
\label{obser-1}
For every graph $G$, $i(G) \le \alpha_c(G) \le \alpha(G)$.
\end{obs}

A graph $G$ is said to be \emph{well-covered\/} if $i(G)=\alpha(G)$. Equivalently, $G$ is well-covered if every maximal independent set of $G$ is a maximum independent set of $G$. The concept of well-covered graphs was introduced by Plummer \cite{Plummer1} in 1970 and extensively studied in many papers. For a~survey of well-covered graphs, see \cite{Plummer}.

A graph $G$ is called an \emph{$\alpha$-excellent graph} if $\alpha_c(G) = \alpha(G)$. Equivalently, $G$ is $\alpha$-excellent if every vertex is contained in some $\alpha$-set of $G$. The $\alpha$-excellent graphs, referred to as $B$-graphs in \cite{Berge1962, Berge1981}, have been studied in \cite{Berge1962, Berge1981, Berge1985, DeLeTo-21, Domke, Fricke, Pushpalatha, Pushpalatha2}.

We are interested in characterizations of graphs $G$ for which $i(G)=\alpha_c(G)$, $i(G)=\alpha(G)$, and $\alpha_c(G)=\alpha(G)$, respectively. We remark that the family of well-covered graphs is properly contained in the family of $\alpha$-excellent graphs, and in the family of graphs $G$ for which  $i(G)=\alpha_c(G)$, respectively. Trees $T$ for which $i(T) = \alpha_c(T)$, and block graphs $G$ such that $\alpha_c(G) = \alpha(G)$ have been characterized in \cite{DeLeTo-21}.
The $\alpha$-excellent trees have also been studied in \cite{Domke, Fricke}.
In this paper, we provide characterizations of $\alpha$-excellent bipartite graphs, unicyclic graphs, simplicial graphs, chordal graphs, block graphs,  and we show that every generalized Petersen graph is $\alpha$-excellent.

\subsection{Notation}

For notation and graph theory terminology, we in general follow~\cite{HaHeHe-22}. Specifically, let $G$ be a graph with vertex set $V(G)$ and edge set $E(G)$, and of order~$n(G) = |V(G)|$ and size $m(G) = |E(G)|$. Two vertices are \emph{neighbors} if they are adjacent.  For a graph $G$ and a vertex $v \in V(G)$, the \emph{open neighborhood} $N_G(v)$ of $v$ is the set of neighbors of $v$ in $G$, and the \emph{closed neighborhood} $N_G[v] = N_G(v) \cup \{v\}$. For a set $S \subseteq V(G)$, the  \emph{open neighborhood} of $S$ is the set $N_G(S) = \bigcup_{v \in S} N_G(v)$, and the \emph{closed neighborhood} $N_G[S] = N_G(S) \cup S$. 

An \emph{end}-\emph{edge} of a graph $G$ is an edge incident with at least one vertex of degree~$1$ in $G$. A \emph{support vertex} is a vertex adjacent to a vertex of degree~$1$, that is, a support vertex is incident with an end-edge. A \emph{strong support vertex} is a vertex with at least two neighbors of degree~$1$. We note that if two end-edges share a vertex, then such a vertex is a strong support vertex. An \emph{interior vertex} of $G$ is a vertex of degree at least~$2$.

A \emph{caterpillar-wheel\/} is a connected unicyclic graph $G$ of girth at least three in which a single cycle $C$ (the \emph{body}) is incident to (or contains) every edge, that is, at least one end of every edge of $G$ belongs to the cycle $C$. An end-edge of a~caterpillar-wheel is its \emph{leg}.

A vertex $v$ of a graph $G$ is a \emph{simplicial vertex} if every two vertices belonging to $N_G(v)$ are adjacent in $G$. Equivalently, a  simplicial vertex is a vertex that appears in exactly one clique of a graph, where a \emph{clique} of a graph $G$ is a maximal complete subgraph of $G$. A clique of a graph $G$ containing at least one simplicial vertex of $G$ is called a \emph{simplex} (plural, \emph{simplexes} or \emph{simplices}) of $G$. We note that if $v$ is a simplicial vertex of $G$, then $G[N_G[v]]$ is the unique simplex of $G$ containing $v$. A~graph $G$ is \emph{simplicial\,} if every vertex of $G$ is a simplicial vertex of $G$ or is adjacent to a simplicial vertex of $G$, that is, a graph $G$ is simplicial if and only if every vertex of $G$ belongs to at least one simplex of $G$.

A graph $G$ is said to be \emph{chordal} (also called \emph{triangulated, rigid-circuit, perfect eliminated,} or \emph{monotone sensitive} in the literature\,) if every cycle of $G$ of length four or more contains a \emph{chord}, i.e., an edge joining two non-consecutive vertices of the cycle. We remark that chordal graphs have  been also characterized in terms of the existence of simplicial vertices \cite{Rose}: A~graph $G$ is  chordal if and only if every induced subgraph of $G$ has a simplicial vertex. Certainly, every induced subgraph of a~chordal graph is chordal.

A \emph{block} in a graph $G$ is a maximal connected subgraph having the property that it contains no cut-vertex of its own. An \emph{end-block} of $G$ is a block containing exactly one cut-vertex, while an \emph{inner-block} of $G$ is a block containing at least two cut-vertices. A graph is a \emph{block graph} if every block of $G$ is a complete graph. It is known that a graph $G$ is a block graph if and only if there exists a unique induced path between any two vertices of $G$.

For a graph $G$ and a family ${\cal H}=\{H_v \colon v \in V(G)\}$ of nonempty graphs indexed by the vertices of $G$, the \emph{corona} $G \circ {\cal H}$ of $G$ and $\cal H$ is the disjoint union of $G$ and $H_v$, $v\in V(G)$, with additional edges joining each vertex $v$ of $G$ to all vertices of $H_v$. Thus, to form the graph $G \circ {\cal H}$ we add a copy of the graph $H_v$ associated with~$v$ for each vertex $v$ of $G$, and join $v$ to all vertices of $H_v$. If all the graphs of the family $\cal H$ are isomorphic to one and the same graph $H$, then we shall write $G \circ H$ instead of $G \circ {\cal H}$. In particular, the corona $G \circ K_1$ is formed from $G$ by adding for each vertex $v$ in $G$ a new vertex $v'$ and the edge $vv'$.

For $k \ge 1$ an integer, we use the standard notation $[k] = \{1,2,\ldots,k\}$ and $[k]_0 = [k] \cup \{0\} = \{0,1,\ldots,k\}$.

\section{Preliminary results}

In this section, we present some preliminary results that we will need when proving our main results. We have the following simple property of simplexes in $\alpha$-excellent graphs.

\begin{prop}
\label{simplexes-are-pairwise-vertex-disjoint}
The simplexes of an $\alpha$-excellent graph are pairwise  vertex-disjoint. In particular, no $\alpha$-excellent graph contains a strong support vertex.
\end{prop}
\begin{proof} Let $G$ be an $\alpha$-excellent graph, and so $\alpha_c(G) = \alpha(G)$. Suppose, to the contrary, that a vertex $v$ of a graph $G$ belongs to two simplexes of $G$, say to $G[N_G[u]]$ and $G[N_G[w]]$. Let $I$ be a maximum independent set of $G$ that contains $v$, and so $\alpha(G) \ge |I| \ge \alpha_c(G) = \alpha(G)$. Consequently, we must have equality throughout this inequality chain, implying that $|I| = \alpha(G)$. However, $(I \setminus \{v\}) \cup \{u,w\}$ is an independent set of $G$ of cardinality~$|I| + 1 = \alpha(G) + 1$, a contradiction.
\end{proof}

Let $G$ be a graph, and let ${\cal H}=\{H_v \colon v \in V(G)\}$ be a family of nonempty graphs indexed by the vertices of $G$. Let $v \in V(G)$, and let $I_v$ be a maximum independent set in $G \circ {\cal H}$ that contains the vertex~$v$. If $H_v$ is not a complete graph, then replacing $v$ in $I_v$ with two non-adjacent vertices in $H_v$ produces an independent set of cardinality greater than~$|I_v|$, implying that $G \circ {\cal H}$ is not $\alpha$-excellent. On the other hand, if $H_v$ is a complete graph for every vertex $v$ of $G$, then $\alpha(G \circ {\cal H}) = |V(G)|$ and every vertex of $G \circ {\cal H}$ belongs to some $\alpha$-set of $G \circ {\cal H}$. This yields the following family of $\alpha$-excellent graphs.

\begin{obs}
\label{corona graph}
Let $G$ be a graph, and let ${\cal H}=\{H_v \colon v \in V(G)\}$ be a family of nonempty graphs indexed by the vertices of $G$. Then the corona $G\circ \cal H$ is an $\alpha$-excellent graph if and only if $\cal H$ consists of complete graphs. In particular, the corona $G \circ K_1$ is an $\alpha$-excellent graph for every graph $G$.
\end{obs}

In \cite{Ravindra}, Ravindra proved that a connected bipartite graph $G$ of order at least 2 is well-covered if and only if $G$ has a perfect matching $M$ and for every $uv \in M$, the induced subgraph $G[N_G(\{u,v\})]$ is a complete bipartite graph. In particular, it was observed in \cite{Ravindra} that a~tree $G$ is well-covered if and only if every interior vertex of $G$ is adjacent to exactly one end vertex of $G$ or, equivalently, if and only if the end-edges of $G$ form a perfect matching of $G$, that is, if and only if $G$ is the corona graph of a tree.
In Proposition \ref{alpha-excellent bipartite graphs} we have an analogue of these results for $\alpha$-excellent bipartite graphs. In fact, this analogue was proved by Berge \cite{Berge1981} many years ago. Our proof given below is a~modification of that given in \cite{Berge1981}. We begin with one definition and two propositions.

A vertex $v$ is a \emph{critical vertex} of $G$ if its removal from $G$ changes the independence number, that is, a vertex $v$ is a critical vertex of $G$ if any of the following equivalent conditions holds: \\[-22pt]
\begin{enumerate}
\item[{\rm (1)}] $\alpha(G-v) \ne \alpha(G)$.
\item[{\rm (2)}] $\alpha(G-v) = \alpha(G) - 1$.
\item[{\rm (3)}] Every maximum independent set of $G$ contains $v$.
\end{enumerate}

\begin{prop}{\rm (\cite{Berge1985})}
\label{Berge: matchings in graphs without critical vertices}
If a graph $G$ contains no critical vertex, then every independent set $I$ of $G$ can be matched into $V(G)\setminus I$.
\end{prop}

\begin{prop}
\label{excellent has no critical vertices}
A connected $\alpha$-excellent graph of order at least~$2$ has no critical vertex, and every independent set $I$ of $G$ can be matched into $V(G)\setminus I$.
\end{prop}
\begin{proof} Let $G$ be a connected $\alpha$-excellent graph of order at least two, and let $v$ be a vertex of $G$. To prove that $v$ is not a critical vertex, it suffices to show that $v$ does not belong to some $\alpha$-set of $G$. Let $u$ be any neighbor of $v$. From the fact that $G$ is an $\alpha$-excellent graph, it follows that $u$ belongs to some $\alpha$-set $I$ of $G$. Since $v$ is adjacent to $u$ and $u\in I$, the vertex $v$ does not belong to $I$, and this proves that $v$ is not a critical vertex of $G$. Consequently, $G$ has no critical vertex and, therefore, by Proposition~\ref{Berge: matchings in graphs without critical vertices}, every independent set $I$ of $G$ can be matched into $V(G)\setminus I$.
\end{proof}

\begin{prop}{\rm (\cite{Berge1981})}
\label{alpha-excellent bipartite graphs}
If $G$ is a connected bipartite graph of order $n\ge 2$, then the following statements are equivalent: \\[-22pt]
\begin{enumerate}
\item[{\rm (a)}] $G$ has a perfect matching.
\item[{\rm (b)}] $\alpha(G) = n/2$. 
\item[{\rm (c)}] $G$ is an $\alpha$-excellent graph.
\end{enumerate}
\end{prop}
\begin{proof}
Assume that $G$ is a connected bipartite graph of order $n\ge 2$, say $G=(A, B, E)$, where $A$~and $B$ be the partite sets of $G$, $|A| \ge |B|\ge 1$, and $E = E(G)$. Since $A$ and $B$ are independent sets and $|A| \ge |B|=n-|A|$, it follows that $\alpha(G)\ge |A|\ge n/2\ge |B|$.

Assume first that $G$ has a perfect matching. Then, certainly, $\alpha(G)\le n/2$ and, therefore, $\alpha(G)= n/2$. Consequently, $|A|=|B|=n/2= \alpha(G)$, and each of the sets $A$ and $B$ is an $\alpha$-set of $G$. Thus every vertex of $G$ belongs to an $\alpha$-set of $G$ and $G$ is an $\alpha$-excellent graph. This proves the implications $(\rm a) \Rightarrow (\rm b)$ and $(\rm b) \Rightarrow (\rm c)$.

Now assume that $G$ is an $\alpha$-excellent graph. We shall prove that $G$ has a perfect matching. By Proposition \ref{excellent has no critical vertices}, $G$ has no critical vertex and every independent set $I$ of $G$ can be matched into $V(G)\setminus I$.
In particular, $A$ can be matched into $B$ and $B$ into $A$. Thus, $|A|=|B|$ and $G$ has a~perfect matching. This proves the implication $(\rm c) \Rightarrow (\rm a)$. \end{proof}

The following simple characterization of $\alpha$-excellent trees follows immediately from Proposition \ref{alpha-excellent bipartite graphs} and it was already proved in \cite{Berge1981}, and independently in \cite{DeLeTo-21, Domke, Fricke}.

\begin{cor} \label{excellent trees}
A tree $T$ of order at least~$2$ is an $\alpha$-excellent graph if and only if $T$ has a~perfect matching. \end{cor}

\section{$\alpha$-Excellent unicyclic graphs}

In this section, we characterize $\alpha$-excellent unicyclic graphs. Our characterization is a counterpart of characterizations of well-covered unicyclic graphs presented in \cite{TV1990,TV1991}.

\begin{prop}
\label{excellent caterpillar-wheels}
Let $G$ be a caterpillar-wheel of girth at least 3. Then $G$ is an $\alpha$-excellent graph if and only if $G$ is a cycle (that is, $G$ is a caterpillar-wheel without any leg) or $G$ is a~caterpillar-wheel with a~perfect matching and with at least two legs.
\end{prop}
\begin{proof}  If $G$ is a cycle, then $G$ is a vertex transitive graph, and, certainly, $G$ is an $\alpha$-excellent graph (as every two vertices of $G$ have the same properties). Thus assume that $G$ is a~caterpillar-wheel with at least two legs and assume that $G$ has a perfect matching. Let $v_1,v_2,\ldots,v_n$ be the consecutive vertices of the only cycle $C$ of $G$, and let $L=\{v_{i_1}\overline{v}_{i_1}, \ldots, v_{i_l} \overline{v}_{i_l}\}$ be the set of legs of $G$, where $1\le {i_1}< \cdots <{i_l}\le n$, and $\overline{L} = \{\overline{v}_{i_1},\ldots, \overline{v}_{i_l}\}$ be the set of leaves of $G$. Let $S$ be the set of support vertices in $G$, that is, $S= \{v_{i_1}, \ldots, v_{i_l}\}$. Let $M$ be the (unique) perfect matching of $G$. Certainly, $L$ is a subset of $M$.

If $M=L$, then $G=C_n\circ K_1$ is a corona graph and it is an $\alpha$-excellent graph, see Observation~\ref{corona graph}. Thus assume that $L$ is a proper subset of $M$. Renaming the vertices on $C$ if necessary, we may assume, without loss of generality, that ${i_1}=1 < \cdots <{i_l}\le n-2$. From the fact that $M$ is a perfect matching of $G$ and $L\varsubsetneq M$, it follows that every component of the subgraph $G'=G-\{v_{i_1}, \overline{v}_{i_1}, \cdots, v_{i_l},\overline{v}_{i_l}\}$ is a path of even order, say $Q_1, \ldots, Q_k$  are the components of $G'$,  where $V(Q_j) = \{u^j_1,u^j_2,\ldots,u^j_{2t_j}\}$ (for some positive integer $t_j$), and the order of vertices of $Q_j$ is inherited from the order of the vertices $v_1,v_2,\ldots,v_n$ on the cycle $C$. It is evident that the perfect matching $M$ is the set
\[
M = L \cup \bigcup_{j=1}^k \{u^j_1u^j_2,\ldots,u^j_{2t_j-1}u^j_{2t_j}\}
\]
and that $\alpha(G)=|M|/2$.

Now we shall prove that $G$ is an $\alpha$-excellent graph. Let $V^j_1$ denote the set $\{u^j_1,\ldots,u^j_{2t_j-1}\}$ of ``odd'' vertices of the path $Q_j$ (with odd subscripts), and, similarly, let $V^j_2$ denote the set $\{u^j_2,\ldots,u^j_{2t_j}\}$ of ``even'' vertices of $Q_j$ (with even subscripts) for $j \in [k]$. It is evident that each of the sets
\[
I_1=\overline{L} \cup \bigcup_{j=1}^k V_1^j \hspace*{0.5cm} \mbox{and} \hspace*{0.5cm} I_2=\overline{L} \cup \bigcup_{j=1}^k V_2^j
\]
is an $\alpha$-set of $G$. Consequently every vertex belonging to $I_1\cup I_2= V(G)\setminus S$ belongs to an $\alpha$-set of $G$. Assume that $v_{i_s}\in S$. It remains to show that $v_{i_s}$ belongs to some $\alpha$-set of $G$. We now consider four cases. In what follows, all subscripts are taken modulo $n$.

\emph{Case 1. $v_{{i_s}-1} \in S$ and $v_{{i_s}+1}\in S$.} In this case, each of the sets $(I_1\setminus \{\overline{v}_{i_s}\})\cup \{v_{i_s}\}$ and $(I_2\setminus \{\overline{v}_{i_s}\})\cup \{v_{i_s}\}$ is an $\alpha$-set of $G$ and $v_{i_s}$ belongs to each of them.

\emph{Case 2. $v_{{i_s}-1} \in S$ and $v_{{i_s}+1}\not\in S$.} In this case, the set
$(I_2\setminus \{\overline{v}_{i_s}\})\cup \{v_{i_s}\}$ is an $\alpha$-set of $G$ that contains $v_{i_s}$.

\emph{Case 3. $v_{{i_s}-1} \not\in S$ and $v_{{i_s}+1}\in S$.} In this case, the set
$(I_1 \setminus \{\overline{v}_{i_s}\})\cup \{v_{i_s}\}$ is an $\alpha$-set of $G$ that contains $v_{i_s}$.

\emph{Case 4. $v_{{i_s}-1} \not\in S$ and $v_{{i_s}+1}\not\in S$.} In this case, without loss generality, assume that $v_{{i_s}-1}$ and $v_{{i_s}+1}$ belong to $Q_k$ and $Q_1$, respectively. In this case, the set
\[
(\overline{L}\setminus \{\overline{v}_{i_s}\}) \cup \{v_{i_s}\} \cup V_1^k \cup \bigcup_{j=1}^{k-1}V_2^j
\]
is an  $\alpha$-set of $G$ that contains $v_{i_s}$. This completes the proof of the first part of the theorem.

Assume now that $G$ is an $\alpha$-excellent caterpillar-wheel which is neither a cycle nor the corona of a cycle. Thus at least one vertex on the cycle $C$ has degree~$2$ in $G$ or at least one vertex on the cycle $C$ is a strong support vertex. By Proposition~\ref{simplexes-are-pairwise-vertex-disjoint}, the $\alpha$-excellent graph $G$ contains no strong support vertex. Hence, at least one vertex on the cycle $C$ is not a support vertex (and has degree~$2$ in $G$). We claim that $G$ has at least two leaves. Suppose, to the contrary, that $G$ has only one leaf. Let $\overline{v}$
be the only leaf of $G$, and let $v$ be the unique neighbor of $\overline{v}$. Adopting our earlier notation, let $C$ denote the (unique) cycle of $G$ and let $C$ have length $n$ where $n \ge 3$. Let $I_v$ and $I_{\overline{v}}$ be a largest independent set of $G$ containing the vertices $v$ and $\overline{v}$, respectively. Since $G$ is an $\alpha$-excellent graph, we have $\alpha(G) = |I_{v}| = |I_{\overline{v}}|$. However, $|I_{\overline{v}}| = 1+\alpha(G-N_G[\overline{v}]) = 1+\alpha(P_{n-1})= 1+\lceil(n-1)/2\rceil$ and $|I_{v}| = 1+\alpha(G-N_G[v]) = 1+\alpha(P_{n-3})= 1+\lceil(n-3)/2\rceil$, and so $|I_{v}| > |I_{\overline{v}}|$, a contradiction. Hence, $G$ has at least two leaves.

Now we shall prove that $G$ has a perfect matching. Let $v$ be a support vertex of $G$, let $\overline{v}$ be a leaf adjacent to $v$, and let $u$ be one of the two neighbors of $v$ on the cycle $C$. Let $G'$ denote the graph $G-uv$. Every independent set of $G$ is an independent set of $G-uv$, and so $\alpha(G-uv)\ge \alpha(G)$. On the other hand, let $J$ be a maximum independent set of $G-uv$. At least one of the sets $J$ or $(J\setminus \{v\})\cup \{\overline{v}\}$ is an independent set of $G$ and, therefore,  $\alpha(G-uv)=|J|= |(J\setminus \{v\})\cup \{\overline{v}\}|\le \alpha(G)$. Consequently, $\alpha(G-uv)= \alpha(G)$. From this and from the fact that $G$ is an  $\alpha$-excellent graph it is evident that $G-uv$ is an $\alpha$-excellent graph. Thus, since $G-uv$ is a tree, Corollary \ref{excellent trees} implies that $G-uv$ has a perfect matching. Therefore, $G$ has a perfect matching.
\end{proof}

Let ${\cal S}$ be the family consisting of all cycles and of all caterpillar-wheels having a perfect matching and at least 2 legs. It follows from the proof of Proposition~\ref{excellent caterpillar-wheels} that a~caterpillar-wheel $G$ which is not a cycle belongs to the family ${\cal S}$ if and only if $G$ has $\ell$ legs, where $\ell \ge 2$, and the distance between any two consecutive legs of $G$ is an odd integer.

Let $\cal H$ be the family of graphs defined recursively as follows: \\ [-22pt]
\begin{enumerate}
\item[{\rm (1)}] The family $\cal H$ contains every graph belonging to the family $\cal S$.
\item[{\rm (2)}] The family $\cal H$ is closed under the operation ${\cal O}_1$ defined below:
\begin{itemize}
\item \textbf{Operation ${\cal O}_1$.} If a graph $H$ belongs to $\cal H$, then add a vertex disjoint copy of a complete graph $K_2$ of order~$2$ to $H$ and add an edge joining a vertex of $H$ with a vertex in the added copy of $K_2$.
\end{itemize}
\end{enumerate}

\begin{prop}
\label{attaching K2 to a graph}
If $H$ is a connected graph of order at least~$2$, and $G$ is a graph obtained from $H$ by applying operation~${\cal O}_1$, then the following properties hold:\\ [-22pt]
\begin{enumerate}
\item[{\rm (1)}] $\alpha(G)=\alpha(H)+1$.
\item[{\rm (2)}] $G$ is $\alpha$-excellent if and only if $H$ is $\alpha$-excellent.
\end{enumerate}
\end{prop}
\begin{proof} Let $v$ and $w$ be the two vertices added to $H$ when constructing $G$ using operation~${\cal O}_1$, and let $u$ be the vertex of $H$ that is joined to exactly one of $v$ and $w$, say $v$. Thus, $uvw$ is a path in $G$, where $w$ has degree~$1$ in $G$ and $v$ has degree~$2$ in $G$. Every $\alpha$-set in $H$ can be extended to an independent set of $G$ by adding to it the vertex~$w$, implying that $\alpha(G) \ge \alpha(H)+1$. On the other hand, every $\alpha$-set in $G$ is an independent set of $G - uv = H \cup K_2$, and therefore $\alpha(G) \le \alpha(G - uv) = \alpha(H)+1$. Consequently, $\alpha(G) = \alpha(H)+1$. This proves part~(1).

To prove part~(2), assume first that $G$ is $\alpha$-excellent. Let $x$ be a vertex of $H$, and let $I_x$ be an $\alpha$-set of $G$ that contains $x$. We note that the set $I_x$ contains exactly one of $v$ and $w$. If $v \in I_x$, then we can replace $v$ in $I_x$ with the vertex~$w$. Hence, we can choose the set $I_x$ to contain the vertex~$w$. The set $I_x \setminus \{w\}$ is an independent set of $H$ of cardinality~$|I_x| - 1 = \alpha(G) - 1 = \alpha(H)$, and is therefore an $\alpha$-set of $H$ that contains the vertex~$x$. Since $x$ is an arbitrary vertex of $H$, the graph $H$ is, therefore, an $\alpha$-excellent graph.

Now assume that $H$ is an $\alpha$-excellent graph. Let $y$ be a vertex of $G$. If $y$ belongs to $V(H)$, and if $I_y$ is an $\alpha$-set of $H$ containing $y$, then $I_y \cup \{w\}$ is an $\alpha$-set of $G$ containing $y$ and $w$. If $t$ is a neighbor of $u$ in $H$ and if $I_t$ is an $\alpha$-set of $H$ containing $t$, then $I_t \cup\{v\}$ is an $\alpha$-set of $G$ containing~$v$. This proves that $G$ is an $\alpha$-excellent graph, and completes the proof of part~(2).
\end{proof}

We define next the opposite operation to Operation ${\cal O}_1$.

\begin{itemize}
\item \textbf{Operation ${\cal O}_2$.} \label{Operation O2} If $G$ is a graph of order at least~$4$ that contains a support vertex of degree~$2$, then remove this support vertex and its (unique) neighbor of degree~$1$ from the graph $G$.
\end{itemize}

We call Operation ${\cal O}_2$ a \emph{plucking operation}. The \emph{plucked graph} of a graph $G$ of order at least~$4$, denoted by $P(G)$, is a graph obtained from $G$ by repeated applications of the plucking operation ${\cal O}_2$ until no further plucking operation is possible. It follows from Corollary \ref{excellent trees} and
from the above definitions that a tree $T$ of order at least~$2$ is an $\alpha$-excellent graph if and only if the plucked graph $P(T)$ is $K_2$. It also follows from these definitions that if a graph $G$ belongs to the family $\cal H$, then $G$ and $P(G)$ are unicyclic graphs. Now from Propositions~\ref{excellent caterpillar-wheels} and~\ref{attaching K2 to a graph} we have the following characterization of the $\alpha$-excellent unicyclic graphs. We omit its straightforward inductive proof.

\begin{thm}
\label{excellent unicyclic graphs}
Let $G$ be a connected unicyclic graph. Then, $G$ is an $\alpha$-excellent graph if and only if one of the following two statements holds:\\ [-22pt]
\begin{enumerate}
\item[{\rm (1)}] The plucked graph $P(G)$ of $G$ is a cycle.
\item[{\rm (2)}] The plucked graph $P(G)$ of $G$ is a caterpillar-wheel having a~perfect matching with $\ell$ legs where $\ell \ge 2$.
\end{enumerate}
\end{thm}

The \emph{subdivision graph} of a graph $G$, denoted $S(G)$, is the graph obtained from $G$ by subdividing every edge of $G$ exactly once. As an immediate consequence of Theorem~\ref{excellent unicyclic graphs}, we have the following result.

\begin{cor}
If $G$ is a connected unicyclic graph, then its subdivision graph $S(G)$ is an $\alpha$-excellent graph.
\end{cor}
\begin{proof} The result is obvious if $G$ is a cycle. If $G$ is a connected unicyclic graph with at least one leg, then $S(G)$ is a unicyclic graph, and its plucked graph $P(S(G))$ is a cycle. Thus, by Theorem~\ref{excellent unicyclic graphs}, $S(G)$ is
an $\alpha$-excellent graph.
\end{proof}

\section{$\alpha$-Excellent simplicial, chordal, and block graphs}

In this section we characterize simplicial graphs, chordal graphs, and block graphs which are $\alpha$-excellent graphs. We begin with the following characterization of well-covered simplicial graphs and well-covered chordal graphs which was proved in \cite{PrToVe-96}.

\begin{prop}{\rm (\cite{PrToVe-96})}
\label{well-covered-simplicial-graphs}
If $G$ is a simplicial graph or a chordal graph, then $G$ is well-covered if and only if every vertex of $G$ belongs to exactly one simplex of $G$.
\end{prop}

The following theorem gives a simple characterization of the $\alpha$-excellent simplicial graphs. This theorem also shows that a simplicial graph is $\alpha$-excellent if and only if it is well-covered.

\begin{thm}
If $G$ is a simplicial graph, then the following statements are equivalent:\\[-22pt]
\begin{enumerate}
\item[{\rm (a)}] $G$ is a well-covered graph.
\item[{\rm (b)}] $G$ is an $\alpha$-excellent graph.
\item[{\rm (c)}] Every vertex of $G$ belongs to exactly one simplex of $G$.
\end{enumerate}
\end{thm}
\begin{proof} The statements (a) and (c) are equivalent, by Proposition \ref{well-covered-simplicial-graphs}. The implication $(\rm a) \Rightarrow (\rm b)$ is obvious, as if $i(G)=\alpha(G)$, then by Observation \ref{obser-1} we have $\alpha_c(G)=\alpha(G)$, and so $G$ is $\alpha$-excellent. Finally, to prove the implication $(\rm b) \Rightarrow (\rm a)$, assume that $G$ is an $\alpha$-excellent graph and suppose, to the contrary, that $G$ is not well-covered. In this case the equivalence of (a) and (c) implies that some vertex of $G$ belongs to at least two simplexes in $G$, which is impossible in an $\alpha$-excellent graph, see Proposition \ref{simplexes-are-pairwise-vertex-disjoint}.
\end{proof}

\begin{figure}[h!] \begin{center} \hspace*{2ex}\includegraphics[scale=0.95]{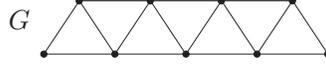}
\caption{An $\alpha$-excellent chordal graph $G$ which is not a well-covered graph}  \label{chordal graph} \end{center}\end{figure}

The graph $G$ shown in Fig. \ref{chordal graph} is an $\alpha$-excellent chordal graph which is not a well-covered graph, as it is easy to observe that $i(G)=2<3 = \alpha_c(G)= \alpha(G)$. Thus, since every well-covered graph is an $\alpha$-excellent graph, the set of well-covered simplicial graphs is  properly contained in the set of $\alpha$-excellent chordal graphs.

For our characterization of $\alpha$-excellent chordal graphs, we need the next two properties of simplicial vertices and the definition of a~successive clique-cover of a graph.

\begin{prop} \label{alpha for G minus neighborhood} If $v$ is a simplicial vertex of a graph $G$, then $\alpha(G-N_G[v])=\alpha(G)-1$.
\end{prop}

\begin{proof} If $I$ is a maximum independent set of $G-N_G[v]$, then $I\cup\{v\}$ is an independent set of $G$ and, therefore, $\alpha(G)\ge |I\cup\{v\}| = |I|+1=\alpha(G-N_G[v])+1$. Now assume that $J$ is a~maximum independent set of $G$. Since $G[N_G[v]]$ is a~complete graph, necessarily $|J\cap N_G[v]|\le 1$ and $J \setminus N_G[v]$ is an independent set of $G-N_G[v]$, implying that $\alpha(G-N_G[v])\ge |J \setminus N_G[v]\ge |J|-1= \alpha(G)-1$. Consequently, $\alpha(G-N_G[v])=\alpha(G)-1$.
\end{proof}

\begin{prop}
\label{alpha excellent for G minus neighborhood}
If $v$ is a simplicial vertex of an $\alpha$-excellent graph $G$, then $G-N_G[v]$ is an $\alpha$-excellent graph.
\end{prop}
\begin{proof} Assume that $v$ is a simplicial vertex of an $\alpha$-excellent graph $G$.
We shall show that $G'= G-N_G[v]$ is an $\alpha$-excellent graph. Let $u$ be an arbitrary  vertex of $G'$. By Proposition~\ref{alpha for G minus neighborhood}, $\alpha(G')=\alpha(G)-1$. Since $G$ is $\alpha$-excellent, in $G$ there is an independent set, say $J_u$, of cardinality $\alpha(G)$ that contains~$u$. The set $J_u\setminus N_G[v]$ is an independent set in $G'$ of cardinality $\alpha(G)-1 = \alpha(G')$ that contains~$u$. Since $u$ is an arbitrary  vertex of $G'$, every vertex of $G'$ belongs to an $\alpha$-set of $G'$, that is, $G'$ is an $\alpha$-excellent graph.
\end{proof}

An ordered sequence $(V_1,\ldots, V_n)$ of subsets of the vertex set $V(G)$ of a graph $G$ is said to be a~\emph{successive clique-cover} of $G$ if the following three properties hold:\\[-22pt]
\begin{enumerate}
\item[{\rm (a)}] The sets $V_1,\ldots, V_n$ form a partition of $V(G)$;
\item[{\rm (b)}] $G[V_i]$ is a simplex in $G-\bigcup_{j<i}V_j$ for $i \in [n]$;
\item[{\rm (c)}] For every $u_i\in V_i$ there are vertices $u_{i+1}\in V_{i+1}, \ldots, u_n\in V_n$ such that $\{u_i, u_{i+1},\ldots, u_n\}$ is an independent set for $i \in [n-1]$.
\end{enumerate}

\begin{prop} \label{a graph with successive clique-cover} If a graph $G$ has a successive clique-cover, then $G$ is an $\alpha$-excellent graph.
\end{prop}

\begin{proof} The result is obvious if $G$ is a complete graph. Thus assume that $G$ is a non-complete graph and  $(V_1,\ldots, V_n)$ is a successive clique-cover of $G$. Let $v_i$ be a simplicial vertex of $G-\bigcup_{j<i}V_j$ belonging to $V_i$ where $i \in [n]$. From the recursive definition of the sequence $(V_1,\ldots, V_n)$ it is evident that the set $I=\{v_1,\ldots,v_n\}$ is independent in $G$. Therefore, $\alpha(G)\ge |I|=n$. On the other hand, let $J$ be a maximum independent set of $G$. Since $G[V_i]$ is a complete graph, we note that $|J\cap V_i| \le 1$ for $i \in [n]$. Thus, since $V_1,\ldots, V_n$ form a partition of $V(G)$, we have $\alpha(G)=|J|= |J\cap V(G)|= |J\cap \bigcup_{i=1}^nV_i| = \sum_{i=1}^n |J\cap V_i|\le n$. Consequently, $\alpha(G)=n$.

To prove that $G$ is an $\alpha$-excellent graph, it suffices to show that every vertex $v$ of $G$ belongs to an independent set of cardinality $n$. This is obvious if $v\in \{v_1,\ldots,v_n\}$ (where, in fact, $v_n$ can be any vertex belonging to $V_n$, as $G[V_n]$ is a complete graph), as $\{v_1,\ldots,v_n\}$ is an independent set of cardinality $n$ in $G$. Thus assume that $v\in V(G)\setminus (\{v_1,\ldots,v_{n-1}\}\cup V_n)$. If $v\in V_1$, then from the property (c) of the successive clique-cover $(V_1,\ldots, V_n)$ of $G$ it follows that there are vertices $u_{2}\in V_{2}, \ldots, u_n\in V_n$ such that $\{v, u_{2},\ldots, u_n\}$ is an independent set of cardinality $n$ that contains~$v$. If $v\in V_i$, where $1<i<n$, then again from the property (c) of the successive clique-cover $(V_1,\ldots, V_n)$ it follows that there are vertices $u_{i+1}\in V_{i+1}, \ldots, u_n\in V_n$ such that $\{v, u_{i+1},\ldots, u_n\}$ is an independent set. In this case $\{v_1,\ldots, v_{i-1}, v, u_{i+1},\ldots, u_n\}$ is an independent set of cardinality~$n$ in $G$ that contains the vertex~$v$. This proves that $G$ is an $\alpha$-excellent graph.
\end{proof}

The next theorem presents a characterization of the $\alpha$-excellent chordal graphs.

\begin{thm} \label{alpha-excellent chordal graph} A chordal graph $G$ is an $\alpha$-excellent graph if and only if $G$ has a successive clique-cover.
\end{thm}
\begin{proof}
If a chordal graph $G$ has a successive clique-cover, then $G$ is an $\alpha$-excellent graph, by Proposition \ref{a graph with successive clique-cover}. We proceed by induction on the order~$n$ to show that an $\alpha$-excellent chordal graph has a successive clique-cover. The result is immediate to check for small $n \le 4$. Let $n \ge 5$ and assume that every $\alpha$-excellent chordal graph of order less than~$n$ has a successive clique-cover. Let $G$ be an $\alpha$-excellent chordal graph of order~$n$.

Let $v_0$ be a simplicial vertex of $G$ and let $G'=G-N_G[v_0]$. The graph $G'$ is a chordal graph of order less than~$n$. In addition, by Proposition~\ref{alpha excellent for G minus neighborhood}, the graph $G'$ is an $\alpha$-excellent graph and, by Proposition \ref{alpha for G minus neighborhood}, $\alpha(G')=\alpha(G)-1$.  Now, by the induction hypothesis, $G'$ has a successive clique-cover, say $(V_1,\ldots,V_{k})$ is a successive clique-cover of $G'$ (where $k=\alpha(G')$), that is, $(V_1,\ldots,V_k)$ is a sequence having the properties (a), (b) and (c) specified in the definition of a successive clique-cover of $G'$. It remains to prove that the sequence $(V_0,V_1,\ldots,V_k)$, where $V_0=N_G[v_0]$, has the desired properties of a successive clique-cover of $G$.

Since $V_1,\ldots,V_k$ is a partition of $V(G')$ and $V_0 = V(G)\setminus V(G')$, we note that $V_0,V_1,\ldots,V_k$ is a partition of $V(G)$. Recall that $v_0$ is a simplicial vertex of $G$ and $V_0 = N_G[v_0]$. We note that $G[V_i]=G'[V_i]$ and $G'[V_i]$ is a simplex in $G' - \bigcup_{j<i}V_j= G - \bigcup_{j<i}V_j$ for $i \in [n]$. These observations imply that $G[V_0]$ is a simplex in $G\,\, (= G- \bigcup_{j<0}V_j)$. Finally, the property (c) of the sequence  $(V_0,V_1,\ldots,V_k)$ is obvious, for if $u_i\in V_i$ and $I_{u_i}$ is a largest independent set containing $u_i$, then, since $G$ is an $\alpha$-excellent graph, we infer that $|I_{u_i}|= k+1$. Moreover, since $G[V_j]$ is a complete graph and $V_0,V_1,\ldots,V_k$ is a partition of $V(G)$, we have $|I_{u_i}\cap V_j|=1$ for $j \in [k]_0$. This proves that $(V_0,V_1,\ldots,V_k)$ a successive clique-cover of $G$.
\end{proof}

It was proved in \cite{TV1990} (see also \cite{PrToVe-96,RanderathVolkmann1994,TV1990Results}) that a~block graph $G$ is well-covered if and only if every vertex of $G$ belongs to exactly one simplex of $G$. The $\alpha$-excellent block graphs have already been characterized in \cite{DeLeTo-21}. Here we present another characterization of such graphs.
In order to state and prove our characterization, we need additional definitions and terminology.

A~family ${\cal P}$ of vertex disjoint blocks of a~block graph $G$ is a \emph{perfect block-cover} of $G$ if every vertex of $G$ belongs to some block in ${\cal P}$. (It is easy to observe that every block graph has at most one perfect block-cover.)   Let $K_n$ be a~complete graph of order $n\ge 1$ with vertex set $\{v_1,\ldots, v_n\}$, and let ${\cal H} =\{H_1,\ldots,H_n\}$ be a~family of $n$ complete graphs. Recall $K_n \circ {\cal H}$  denotes a graph obtained from the disjoint union of the graphs $K_n, H_{1},\ldots, H_{n}$ by joining the vertex $v_i$ of ${K_n}$ to every vertex of $H_{i}$ for $i \in [n]$. The graph $K_n\circ {\cal H}$ is said to be the {\em general corona} of $K_n$ and ${\cal H}$, and its complete subgraph induced by the vertices $v_1,\ldots, v_n$ is called the {\em body} of $K_n\circ {\cal H}$.

Let ${\cal F}$ be the family of graphs defined in~\cite{DeLeTo-21} and such that: \\ [-22pt]
\begin{enumerate}
\item[{\rm (1)}] contains every complete graph of order at least 2; and
\item[{\rm (2)}] is closed under attaching general coronas of complete graphs
\end{enumerate}

\begin{figure}[h!] \begin{center} \hspace*{2ex}\includegraphics[scale=0.9]{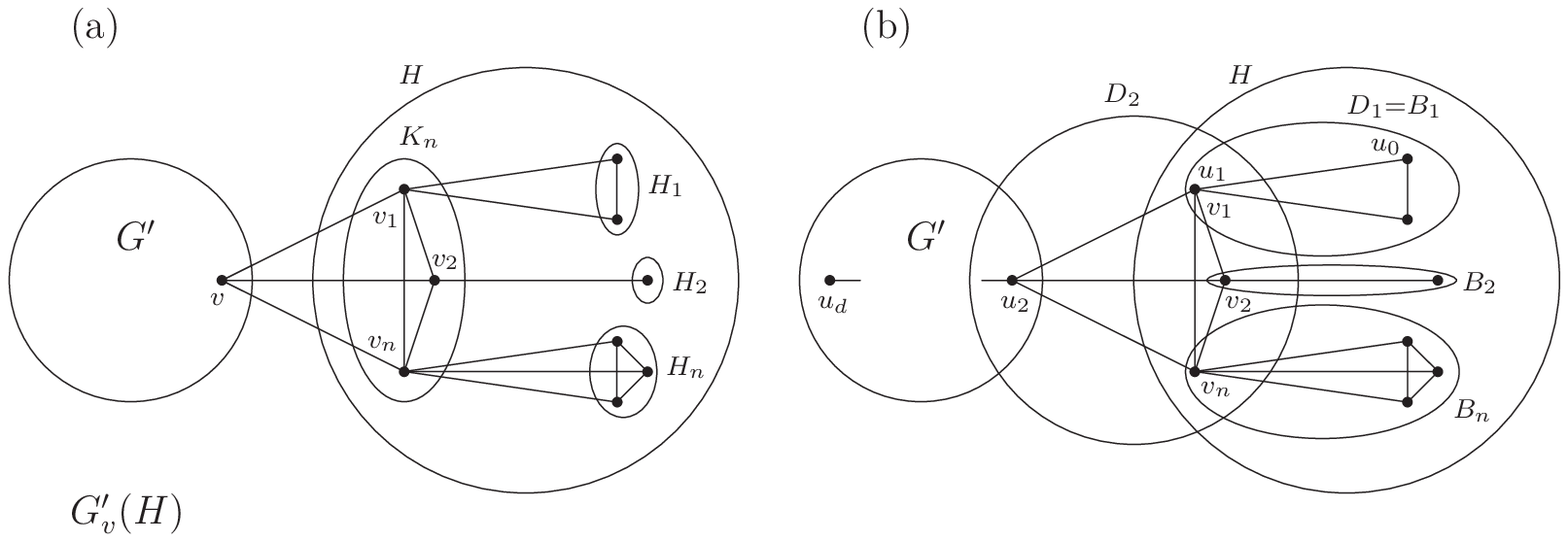}
\caption{} \label{fig-block-graphs} \end{center}\end{figure}

Thus, if a graph $G'$ belongs to the family ${\cal F}$  and $H=K_n\circ {\cal H}$ is a general corona of a complete graph $K_n$ and a family $\cal H$ of complete graphs,  then to ${\cal F}$ belongs every graph $G'_v(H)$ obtained from the disjoint union $G'\cup H$ by adding the edges joining the body of $H$ with the vertex $v$ of $G'$, see Fig. \ref{fig-block-graphs}\,(a). It is clear from this definition that every graph belonging to the family ${\cal F}$ is a~block graph.

We are now in a position to state and prove a characterization of the $\alpha$-excellent block graphs in terms of perfect block-covers.

\begin{thm}
\label{m1}
Let $G$ be a connected block graph of order at least $2$. Then the following statements are equivalent:\\[-22pt]
\begin{enumerate}
\item[{\rm (a)}] $G$ is an $\alpha$-excellent graph.
\item[{\rm (b)}] $G\in \mathcal{F}$.
\item[{\rm (c)}] $G$ has a perfect block-cover.
\end{enumerate}
\end{thm}
\begin{proof} The equivalence between $(\rm a)$ and $(\rm b)$ was proved in \cite{DeLeTo-21}. By making a little modification of that proof, we shall prove the equivalence of the statements $(\rm b)$ and $(\rm c)$. Assume first that $G\in {\cal F}$. Thus, $G$ can be obtained from a~sequence $G_1,\ldots, G_m$ of graphs belonging to ${\cal F}$, where $G_1$ is a~complete graph of order at least~$2$ and $G=G_m$, and, if $m\ge 2$, $G_{i+1}$ can be obtained from $G_i$ by attaching a general corona (of a~complete graph) for $i \in [m-1]$. By induction on the number $m$ we shall prove that $G$ has a~perfect block-cover.

If $m=1$, then $G$ is a complete graph of order at least~$2$ and, certainly, the only block of $G$ forms a~perfect block-cover of $G$. This establishes the base case of the induction. Assume, then, that the result holds for all graphs belonging to ${\cal F}$ that can be constructed from a sequence of fewer than $m$ graphs, where $m\ge 2$. Let $G$ be obtained from a sequence $G_1,\ldots,G_{m-1}$ by attaching a general corona $H=K_n\circ\{H_1,\ldots,H_n\}$ to a vertex $v$ of $G'=G_{m-1}$, that is, assume that $G= G'_v(H)$. By the induction hypothesis, $G'$ has a perfect block-cover, say ${\cal P}'$. Now, it is obvious that ${\cal P}'\cup \{B_1, \ldots,B_n\}$ is a~perfect block-cover of $G$, where $B_1, \ldots,B_n$ are the blocks of $G$ (and of $H$) induced by the vertex-sets $V(H_1)\cup \{v_1\}, \ldots, V(H_n)\cup \{v_n\}$, respectively, where $v_1,\ldots, v_n$ are the vertices of the body of $H$. This proves the implication $(\rm b) \Rightarrow (\rm c)$.

Finally, to prove the implication $(\rm c) \Rightarrow (\rm b)$, assume that $G$ is a block graph of order at least~$2$ and $G$ has a~perfect block-cover, say  ${\cal P}=\{B_1,\ldots,B_n\}$. We proceed by induction on the number $n$ of blocks in ${\cal P}$. If $n=1$, then $G=B_1$ is a complete graph and, certainly, $G\in {\cal F}$. Assume, then, that $n\ge 2$ and that the statement holds for all block graphs having a perfect block-cover of cardinality less than~$n$. From the fact that the blocks $B_1,\ldots,B_n$ are disjoint and $n\ge 2$, it follows that the diameter $d$ of $G$ is greater than 2. Let $P \colon u_0u_1 \ldots u_d$ be a longest path without chords in $G$. Let $D_1$ and $D_2$ be the blocks in $G$ that contain the edges $u_0u_1$ and $u_1u_2$, respectively. From the choice of $P$ as the longest chordless path in $G$, it follows that $u_0$ is a simplicial vertex in $G$. Therefore, $D_1$ is the only block in $G$ that contains $u_0$. The block $D_1$, as does every other simplex of $G$, belongs to ${\cal P}$, and, consequently, every two simplexes of $G$ are pairwise vertex-disjoint. This also means that $D_1$ is the only simplex of $G$ that contains $u_1$.  We now consider two cases depending on the order of $D_2$.

\medskip
\emph{Case 1. $D_2$ is of order two.} In this case, ${\cal P}'={\cal P}-\{D_1\}$ is a~perfect block-cover of the subgraph $G'=G-V(D_1)$ of $G$. Thus, by the induction hypothesis, $G'\in {\cal F}$. Consequently, $G$ belongs to the family $\cal F$ as $G$ can be obtained from $G'$ by attaching the general corona $K_1\circ\{D_1-u_1\}=G[\{u_1\}] \circ\{D_1-u_1\}$ to the vertex $u_2$ in $G'$.

\emph{Case 2. $D_2$ is of order at least three.} Assume that $u_2, v_1=u_1,v_2,\ldots, v_l$ are the vertices of $D_2$, where $l\ge 2$. Since $D_1$ and $D_2$ share a vertex, the block $D_2$ is not a simplex. Thus every vertex of $D_2$ is a cut-vertex in $G$. From this and from the choice of $P$ as the longest chordless path in $G$, it follows that each vertex $v_i$ belongs to some end-block of $G$. In addition, since every end-block of $G$ belongs to ${\cal P}$ and blocks in ${\cal P}$ are disjoint, every vertex $v_i$ of $D_2$ belongs to exactly one end-block in $G$, say $B_i$ is an end-block of $G$ that contains $v_i$. Certainly, $B_1=D_1$, and the subgraph $H$ of $G$ induced by the vertices belonging to the blocks $B_1, B_2,\ldots,B_l$ is a general corona, $H= K_l \circ {\cal H}$, where $K_l=G[\{v_1,\ldots,v_l\}]$ and ${\cal H}=\{B_1-v_1, \ldots, B_l-v_l\}$. Now, ${\cal P} \setminus \{B_1, \ldots,B_l\}$ is a~perfect block-cover of the subgraph $G'=G-V(H)$ of $G$. Thus by the induction hypothesis, $G'\in {\cal F}$. Consequently, $G$ belongs to the family $\cal F$ as $G$ can be obtained from $G'$ by attaching the general corona $H$ to the vertex $u_2$ in $G'$, see Fig. \ref{fig-block-graphs}\,(b). This completes the proof of Theorem~\ref{m1}.
\end{proof}

\section{$\alpha$-Excellent generalized Petersen graphs}

Let $n$ and $k$ be positive integers with $n\ge 3$ and $k \in [n-1]$. The
\emph{generalized Petersen graph} $P_{n,k}$ is defined in the following way. It
has $2n$ vertices, say $v_0, v_1,\ldots,v_{n-1}$, $u_0, u_1,\ldots,u_{n-1}$, and
edges $v_iv_{i+1}$, $v_iu_i$, and $u_iu_{i+k}$ for all $i$ satisfying $i \in [n-1]_0$ with all subscripts taken modulo $n$. By construction each vertex $v_i$ is of degree~$3$ in $P_{n,k}$. Similarly, each $u_i$ is a vertex of degree~$3$ if $k \ne n/2$ but its degree is~$2$ if $k=n/2$. Since $P_{n,k}$ is isomorphic to $P_{n,n-k}$, we may therefore always assume that $k\le \lfloor n/2\rfloor$. A simple analysis shows that $P_{n,k}$ is a graph of girth~$3$ if and only if $n=3k$ for $k\ge 1$. Analogously, $P_{n,k}$ is a graph of girth~$4$ if and only if $k=1$ and $n\ge 4$, $k=2$ and $n=4$ or $k\ge 1$ and $n=4k$.

\begin{figure}[h!] \begin{center} \hspace*{2ex}\includegraphics[scale=0.9]{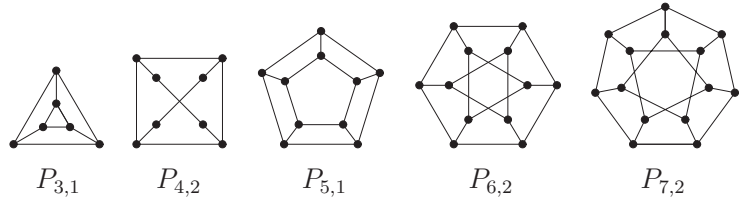}
\caption{The well-covered generalized Petersen graphs} \label{Petersen graphs}\end{center}\end{figure}

It was shown in \cite{Topp1995} that $P_{3,1}$, $P_{4,2}$, $P_{5,1}$, $P_{6,2}$, and $P_{7,2}$ (shown in Fig. \ref{Petersen graphs}) are the only well-covered generalized Petersen graphs. In contrast, all generalized Petersen graphs are $\alpha$-excellent graphs.

\begin{prop}
Every generalized Petersen graph is an $\alpha$-excellent graph.
\end{prop}
\begin{proof}
Let $G$ be the generalized Petersen graph $P_{n,k}$ with vertex set $V(G) = V \cup U$, where $V = \{v_0, v_1,\ldots,v_{n-1}\}$, $U=\{u_0, u_1, \ldots, u_{n-1}\}$, and edges $v_iv_{i+1}$, $v_iu_i$, and $u_iu_{i+k}$ for all $i$ satisfying $i \in [n-1]_0$ with all subscripts taken modulo $n$. Let $I_0$ be an $\alpha$-set of $G$. If $I_0 \cap V = \emptyset$, then $I_0 = U$, which is impossible as $U$ is not an independent set in $G$. Hence, $I_0 \cap V \ne \emptyset$. If $I_0 \cap U = \emptyset$, then $I_0 = V$, which is impossible as $V$ is not an independent set in $G$. Hence, $I_0 \cap U \ne \emptyset$.
Now, if $I_l=\{v_{i+l}\colon v_i\in I_0\cap V\}\cup \{u_{j+l} \colon u_j\in I_0\cap U\}$, where $l \in [k-1]$ and all subscripts are taken modulo $n$, then
each of the sets $I_0, I_1, \ldots, I_{k-1}$ is an $\alpha$-set of $G$, and, certainly,  each vertex of $G$ belongs to at least one of the sets $I_0, I_1, \ldots, I_{k-1}$. This proves that $G$ is an $\alpha$-excellent graph.
\end{proof}

\section{Concluding remarks}

\begin{figure}[h!] \begin{center} \hspace*{2ex}\includegraphics[scale=0.9]{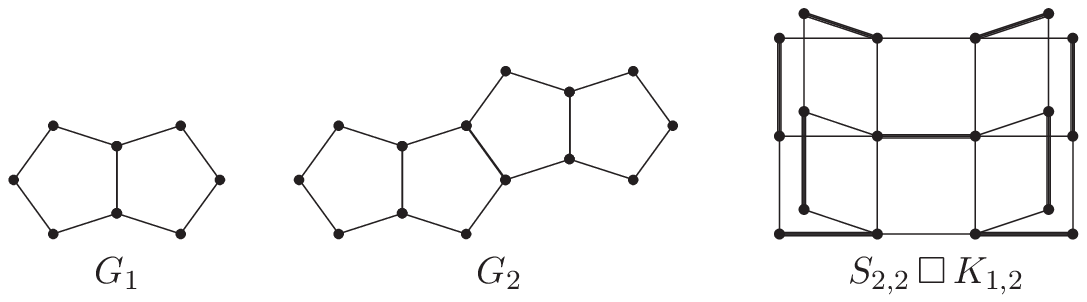}
\caption{Graphs $G_1$, $G_2$, and $S_{2,2} \cp  K_{1,2}$}\label{Graphs C(5)-trees} \end{center}\end{figure}

A graph $G$ is a $C_{(n)}$-\emph{tree} if it can be constructed from a cycle of length $n$ by a finite number of applications of the following operation: add a new cycle of length $n$ and identify an edge of this cycle with an edge of the existing graph. Note that every $2$-tree of order at least~$3$ is a $C_{(3)}$-tree and vice versa. It is known (see \cite{Topp1995}) that a $C_{(n)}$-tree $G$ is a well-covered graph if and only if $G \in \{C_3, C_4, C_5, C_7\}$ or $G$ is a $C_{(3)}$-tree in which every vertex belongs to exactly one end cycle of $G$ (that is, to exactly one simplex of $G$). The way in which each $C_{(n)}$-\emph{tree} is constructed implies that if $n$ is even, then every $C_{(n)}$-\emph{tree} is a~bipartite graph and has a~perfect matching. From this and from Proposition \ref{alpha-excellent bipartite graphs} we conclude that every $C_{(n)}$-tree is an $\alpha$-excellent graph if $n$ is even. However, we do not know which $C_{(n)}$-trees are $\alpha$-excellent if $n$ is odd. The graphs $G_1$ and $G_2$ in Fig. \ref{Graphs C(5)-trees} are $C_{(5)}$-trees, and only $G_2$ is an $\alpha$-excellent graph.

It would also be interesting to know which Cartesian products of graphs are $\alpha$-excellent graphs. The problem is not clear even for bipartite graphs. It is obvious that the Cartesian product $G \cp H$ is a~bipartite graph if and only if both $G$ and $H$ are bipartite and, as we know, a bipartite graph is an $\alpha$-excellent graph if and only if it has a perfect matching. Thus, if $G$ and $H$ are bipartite graphs and at least one of them is an $\alpha$-excellent graph, then the Cartesian product $G \cp H$ is an $\alpha$-excellent graph. The opposite implication is not true, as, for example, the Cartesian product $S_{2,2} \cp  K_{1,2}$ in Fig. \ref{Graphs C(5)-trees} is an $\alpha$-excellent graph but neither the double star $S_{2,2}$ nor the star $K_{1,2}$ is an $\alpha$-excellent graph. 

We close this paper with the following two open problems that we have yet to settle.\\ [-22pt]
\begin{enumerate}
\item[(1)] Characterize the $\alpha$-excellent $C_{(n)}$-trees for odd $n\ge 3$.
\item[(2)] Study $\alpha$-excellent Cartesian products of graphs.
\end{enumerate}

\end{document}